\documentclass[11pt]{amsart}

\usepackage{epigamath-voisin}

\usepackage[english]{babel}

\numberwithin{equation}{section}

\usepackage[shortlabels]{enumitem}
\setlist[enumerate,1]{label={\rm(\arabic*)}, ref={\rm\arabic*}}

\usepackage{amsmath,amssymb,epsf,graphics}
\usepackage{amsfonts,amsxtra,amscd,verbatim,eucal}
\usepackage{mathtools}
\usepackage[all]{xy}

\usepackage{extarrows} 

\newtheorem{thm}{Theorem}[section]
\newtheorem{lemma}[thm]{Lemma}

\theoremstyle{remark}
\newtheorem{rem}[thm]{Remark}

\theoremstyle{definition}
\newtheorem{defn}[thm]{Definition}

\newcommand*{\MyDef}{\,\mathrm{def}\,}
\newcommand*{\eqdef}{\mathop{\overset{\MyDef}{\scalebox{1.8}[1]{=}}}}

\DeclareMathOperator{\vol}{vol}
\renewcommand{\Im}{\operatorname{Im}}
\DeclareMathOperator{\Sym}{Sym}
\DeclareMathOperator{\rk}{rk}
\DeclareMathOperator{\Weil}{Weil}
\DeclareMathOperator{\NS}{NS}

\newcommand{\supth}[1]{\ensuremath{#1^{\mathrm{th}}}}

\YearArticle{2024} \EpigaArticleNr{14} \ReceivedOn{January 9, 2023}
\InFinalFormOn{August 28, 2023}
\AcceptedOn{October 17, 2023}

\title{Coholomogy class of complex approximable algebras}
\titlemark{Coholomogy class of complex approximable algebras}

\author{Catriona Maclean}
\address{Institut Fourier, Universit\'e Grenoble Alpes, CS 40700, 
  38058 Grenoble cedex 09, France}
\email{catriona.maclean@univ-grenoble-alpes.fr} 

\authormark{C.~Maclean}

\AbstractInEnglish{In his work on arithmetic Fujita approximation, 
  Huayi Chen introduces the notion of an approximable graded algebra,
  which he uses to prove a Fujita-type theorem in the arithmetic
  setting, and asks if any such algebra is the graded ring of a big
  line bundle on a projective variety. This was proved to be false by
  the author in a previous work. A subsequent paper showed that, whilst the
  approximable algebra is not necessarily a subalgebra of the algebra
  of graded sections of a big line bundle, it is a full-dimensional
  subalgebra of the algebra of sections of an infinite Weil
  divisor. This paper also proved that over $\mathbb{C}$, if an
  infinite Weil divisor has a finite numerical class, then its section
  ring is approximable, and the same is true for any full-dimensional
  subalgebra.

  This note proves that the converse is true: any
  approximable algebra over $\mathbb{C}$ is associated to an infinite
  Weil divisor whose numerical class converges.}

\MSCclass{14E99}

\KeyWords{Birational geometry, valuation theory, conplex geometry}

\acknowledgement{The author of this paper was supported by the ERC
  grant ALKAGE.}

\begin{document}

%%%%%%%%%%%%%%%%%%%%%%%%%%%%%%%
% Title page
%%%%%%%%%%%%%%%%%%%%%%%%%%%%%%%

\maketitle

\begin{prelims}

\DisplayAbstractInEnglish

\bigskip

\DisplayKeyWords

\medskip

\DisplayMSCclass

\end{prelims}

%%%%%%%%%%%%%%%%%%%%%
% Table of Contents
%%%%%%%%%%%%%%%%%%%%%

\newpage

\setcounter{tocdepth}{1}

\tableofcontents

%%%%%%%%%%%%%%%%%%%%%
% Content begins here
%%%%%%%%%%%%%%%%%%%%%

\section{Introduction}
The Fujita approximation theorem, \cite{fuj}, is an important result
in algebraic geometry. It states that whilst the section ring
associated to a big line bundle $L$ on an algebraic variety $X$
\[ R(L)\eqdef\oplus_m H^0(mL, X)\]
is typically not a finitely generated algebra, it can be approximated
arbitrarily well by finitely generated algebras. More precisely, we have the following. 

\begin{thm}[Fujita]
Let $X$ be an algebraic variety, and let $L$ be a big line bundle on
$X$.  For any $\epsilon>0$ there exist a birational modification
\[\pi\colon \hat{X}\longrightarrow X\]
and a decomposition of\, $\mathbb{Q}$-divisors $\pi^* (L)= A+E$ such that
\begin{itemize}
\item $A$ is ample and $E$ is effective,  
\item $\vol(A)\geq (1-\epsilon)\vol(L)$.
\end{itemize}
\end{thm}

In \cite{LM}, Lazarsfeld and Musta\c{t}\u{a} use the Newton--Okounkov body
associated to $A$ to give a simple proof of Fujita
approximation.

In \cite{volumes_chen},  Chen uses
Lazarsfeld and Musta\c{t}\u{a}'s work on Fujita approximation to prove a
Fujita-type approximation theorem in the arithmetic setting.  In the
course of this work, he defines the notion of approximable graded
algebras, which are exactly those algebras for which a Fujita-type
approximation theorem holds.

\begin{defn}
An integral graded algebra ${\bf B}=\oplus_m B_m$ with $B_0=k$ a field
is approximable if and only if the following conditions are satisfied: 
\begin{enumerate}
\item All the graded pieces $B_m$ are finite-dimensional over $k$.
\item For all sufficiently large $m$ the space $B_m$ is non-empty. 
\item  
For any $\epsilon$ there exists a $p_0$ such that for all $p\geq p_0$
we have that
\[ \liminf_{n\rightarrow \infty} \frac{ \dim(\Im(S^nB_p\rightarrow B_{np}))}{\dim(B_{np})}> (1-\epsilon). \]
\end{enumerate}
\end{defn}

In his paper \cite{volumes_chen}, Chen asks whether any graded
approximable algebra is in fact a subalgebra of the algebra of
sections of a big line bundle.  A counter-example to this is given in
\cite{maclean}, where the graded approximable algebra is equal to the
section ring of an {\it infinite} divisor.\footnote{Infinite in this
context means an infinite formal sum of Weil divisors with real
coefficients $\sum_i a_i D_i$.} The subsequent paper \cite{maclean2}
shows that any approximable algebra is indeed a subalgebra of the
section ring of such a divisor.  Moreover, it is established that if
$X$ is a smooth complex algebraic variety of dimension $d$ and
$D=\sum_{i=1}^\infty a_i D_i$ is an infinite Weil divisor on $X$ such
that the sum of divisor classes $ \sum_{i=1}^\infty a_i[D_i ]$
converges in the N\'eron--Severi space $\NS(X)$ to a finite real big
cohomology class $[D]$, then the algebra $\oplus_m H^0(mD)$ is
approximable.

The main result of this note proves that the converse is true.

\begin{thm}\label{main}
Let ${\bf B}$ be an approximable algebra over $\mathbb{C}$, and let
$X({\bf B})$ and $D({\bf B})=\sum_{i=1}^\infty a_iD_i$ be the smooth
complex variety and infinite Weil divisor constructed in
\cite{maclean2} such that ${\bf B}$ is a full-dimensional subalgebra
of the section ring of\, $D({\bf B})$. The sum of cohomology classes
$\sum_{i=1}^\infty a_i [D_i]$ in $\NS(X)$, the N\'eron--Severi space
of\,~$X$, is then a convergent series.
\end{thm}

In Section 2, we set notation and recall such results from
\cite{maclean2} as will be necessary. Section 3 contains the proof of
Theorem~\ref{main}.

\section{Notation and two preliminary results}
In this section, we  fix some notation and recall some essential preliminary results.

In this paper, we consider a graded approximable algebra ${\bf
  B}=\oplus_m B_m$ such that $B_0=\mathbb{C}$.  For any natural
numbers $k$ and $n$ we denote by $\Sym^n(B_k)$ the $\supth{n}$
symmetric power of the vector space $B_k$ and by $S^n(B_k)$ the image
of $\Sym^n(B_k)$ in $B_{nk}$.

We  say that a graded, not \textit{a priori} approximable, algebra
$\oplus_m B_m$ {\it is of dimension $d$ and has volume $v$} if the
sequence
\[ \lim_{n\rightarrow \infty}
\left( \frac{\rk (B_n)}{ n^{d}/ d!}\right)\]
converges to the strictly positive real number $v$. It is proved by Chen in \cite{volumes_chen} that such numbers exist for any approximable algebra ${\bf B}$.

We now recall  the construction from \cite{maclean2} of the smooth
variety $X({\bf B})$ and the infinite Weil divisor $D({\bf B})$ such
that ${\bf B}$ is a subalgebra of the section ring of $D({\bf B})$.

\subsection{Construction of $\boldsymbol{X({\bf B})}$ and $\boldsymbol{D({\bf B})}$}
The variety $X({\bf B})$ is defined up to birational equivalence using
the homogeneous field of fractions of ${\bf B}=\oplus_m B_m$, which we
now define.

\begin{defn}
Let ${\bf B}=\oplus_m B_m$ be a graded algebra over $\mathbb{C}$. Then we define its homogeneous 
fraction field by
\[ K^{\hom}({\bf B})= \left\{ \frac{b_1}{b_2}\mid \exists m \mbox{ such that } b_1, b_2\in B_m, b_2\neq 0\right\}/ \sim \]
where $\sim$ is the equivalence relation 
\[ \frac{b_1}{b_2}\sim \frac{c_1}{c_2} \Longleftrightarrow b_1 c_2= c_1 b_2.\]
\end{defn}

Note that $K^{\hom}({\bf B})$ is a field extension of
$\mathbb{C}$.

Choose $n$ large enough that $B_n$ and $B_{n+1}$
are both non-trivial. Choose $f_1\in B_n$ and $f_2\in B_{n+1}$. For
any $m$ we can then identify $B_m$ with a subspace of $K^{\hom}({\bf B})$ via the inclusion
\[i_m\colon B_m\xhookrightarrow{\hphantom{aaa}} K^{\hom}({\bf B}),\quad 
b_m\longmapsto \frac{b_m f_1^m}{f_2^m}.\]
(Note that this inclusion depends on the choice of $f_1$ and $f_2$). 
Throughout what follows, we  assume that we have fixed $f_1\in B_n$ and $f_2\in B_{n+1}$ for some choice of $n$ and that
there is therefore for every $m$ a fixed inclusion $i_m\colon B_m\hookrightarrow K^{\hom}({\bf B})$. These inclusions satisfy
$i_m(b_m)i_r(b_r)= i_{m+r} (b_mb_r)$ for any $b_m\in B_m$ and $b_r\in B_r$.

  This field is proved in \cite{maclean2} to be finitely generated,
enabling the following definition of $X({\bf B})$.

  \begin{defn}
The variety $X({\bf B})$ is a smooth projective complex variety such that
$K\left(X({\bf B})\right)= K^{\hom}({\bf B})$.
\end{defn}

  \begin{rem}
The variety $X({\bf B})$ is defined only up to birational equivalence. 
It can be chosen 
smooth by Hironaka resolution.
\end{rem}
  
We now recall the definition of the infinite Weil divisor $D({\bf B})$,
which is constructed as the limit of the sequence of divisors $D_m/m$,
where the divisors $D_m$ are poles of the rational functions $b_m\in
B_m$. More precisely, for any $b_m\in B_m$ let
$(b_m)_X=\sum_{i=1}^\infty a_i D_i$ be the principal divisor\footnote{
Here, the $D_i$ are prime divisors, and the sum is finite.} on $X({\bf
  B})$ cut out by the rational function $b_m$. We let
$(b_m)_X^-=\sum_{i|c_i<0}-c_iD_i$ be the poles divisor
of $b_m$ and let $(b_m)_X^+=\sum_{i|c_i>0} c_i D_i$ be the zeros
divisor of $b_m$, so that
\[ (b_m)_X= (b_m)_X^+-(b_m)_X^-.\] We can now define $D_m$.

\begin{defn}
For any $m$ such that $B_m$ is non-empty, we define the effective
divisor $D_m$ on $X({\bf B})$ by
\[ D_m= \sup_{b_m\in B_m}\left((b_m)_X^-\right),\]
where the supremum is taken with respect to the natural partial order on ${\Weil}(X({\bf B}))$.
\end{defn}

It is proved in \cite{maclean2} that this supremum is a maximum and
moreover that for every $m$ such that $B_m$ is non-empty, there exists a
$b_m$ such that $D_m$ is the poles divisor of $i_m(b_m)$.  We can now
define $D({\bf B})$.

\begin{defn}
We set $D({\bf B})= \lim_{m\rightarrow \infty}
\left(\frac{D_m}{m}\right)$.\footnote{It is proved in \cite{maclean2}
that this limit exists and that for any $k$ the divisor $\lfloor kD\rfloor$
is finite.}
\end{defn}
\noindent

In the next section, we prove that the infinite divisor $D(\mathbf{B})$
has a finite cohomology class.

\section{Finiteness of the cohomology class \texorpdfstring{$\boldsymbol{[D({\bf B})]}$}{[D(B)]}}
We now give a proof of Theorem~\ref{main}. Note that if ${\bf B}_q=
\oplus_n B_{nq}$, then $D({\bf B}_q)= q D({\bf B})$, so without loss of
generality, we may assume that $B_1$ is non-zero. Our proof depends on the following observation.

\begin{lemma}\label{cone}
Let $X$ be a smooth complex variety, and let $H$ be an ample divisor on
$X$. Equip $\NS(X)$, the N\'eron--Severi group of\, $X$, with a norm denoted by 
$| \cdot |$.  Then there is  a constant $C> 0$ such that for any
pseudo-effective divisor $E$ on $X$, we have that
\[ E\cdot H^{d-1} > C | [E] |.\]
\end{lemma}

\begin{proof}
Suppose not. Then there exists a sequence of pseudo-effective divisor
classes $[E_n]$ such that $|[E_n]|=1$ for all $n$ and $E_n\cdot
H^{d-1} \rightarrow 0$. Passing to a convergent subsequence, we may
assume that in the N\'eron--Severi space, $[E_n]$ converges to a non-zero
pseudo-effective cohomology class $[E]$ such that $H^{d-1}\cdot E =0$,
but this is impossible. This completes the proof of Lemma
\ref{cone}.
\end{proof}

This lemma enables us to give a numerical criterion for convergence of
the sequence $[D_m/m]$.

\begin{lemma}\label{num_conv}
Let ${\bf B}$ be an approximable algebra over $\mathbb{C}$, let
$X({\bf B})$ be the associated variety, and consider the divisors
$D_m/m$ defined above, converging to an infinite Weil divisor $D({\bf
  B})$. Let $\NS(X)$ be the N\'eron--Severi space of\,~$X$.  Let $H$ be an
ample divisor on $X$.

The sequence $[D_m/m]$ converges in $\NS(X)$ if and only if the
numerical sequence $(D_m\cdot H^{d-1})/m$ converges.
\end{lemma}

\begin{proof}
The fact that convergence of $[D_m/m]$ implies convergence of
$(D_m\cdot H^{d-1})/m$ is immediate. We assume now that the sequence
$(D_m\cdot H^{d-1})/m$ converges, and we will show that the sequence of
cohomology classes $[D_m/m] $ also converges. Throughout what follows,
we write $[D_1]\geq [ D_2]$ if and only if $[D_1-D_2]$ is the
cohomology class of a {\it pseudo-effective} divisor.

We first show  that if the sequence $[D_{m!}/(m!)]$ converges to a
limit $[D]$, then so does the sequence $[D_{m}/m]$. For any integer
$m_1$ and any $\epsilon > 0$, there exists an integer $M_1$ such that
if $m> M_1$, then
\begin{equation}\label{equ1} [D_m/m]\geq (1-\epsilon) [D_{m_1!}/(m_1!)].\end{equation} We also have that 
\begin{equation}\label{equ2} [D_m/m] \leq [D_{m!}/(m!)]\leq [D]\end{equation} since the sequence $[D_{m!}/(m!)]$ is increasing in $m$. (The last equality is immediate because we have defined the partial order $\leq $ in terms of
pseudo-effective divisors rather than effective divisors).

It follows from (\ref{equ1}) that for all $m> M_1$ we have that
\begin{equation*}
  [D- D_m/m] \leq [(D-(1-\epsilon) D_{m_1!}/(m_1!)]
\end{equation*}
and hence
\[H^{d-1}\cdot (D- D_m/m) \leq [H]^{d-1}\cdot [D-(1-\epsilon) D_{m_1!}/(m_1!)]).\] 
Since the right-hand side can be made arbitrarily small by an
appropriate choice of $\epsilon$ and $m_1$, we deduce that
\begin{equation*}
[H]^{d-1}\cdot ([D- D_m/m])\longrightarrow_{m\rightarrow 0} 0.
\end{equation*}
Since $[D-D_m/m]$ is pseudo-effective by Equation (\ref{equ2}),  it follows from Lemma~\ref{cone} that
\begin{equation*}
  [D_m/m]\longrightarrow_{m\rightarrow \infty} [D].
  \end{equation*}

It remains only to show that the sequence $D_{m!}/m!$ is
convergent. We note that this sequence of divisors is monotone
increasing. In particular, if we set
 \[ R_{m}= D_{m!}/(m!) - D_{(m-1)!}/(m-1)!,\]
then $R_m$ is effective and hence pseudo-effective for all $m$. We
have assumed that the sequence $\sum_{m=1}^\infty H^{d-1}\cdot R_m$ is
convergent. By Lemma~\ref{cone}, the series $\sum_{m=1}^\infty [R_m]$
is also convergent.  This completes the proof of Lemma~\ref{num_conv}.
\end{proof}

 A final lemma will be necessary before completing the proof of Theorem~\ref{main}.

\begin{lemma}\label{division}
There are constants $p$, $k$ and $N$ such that for any $m> N$
divisible by $p$ and any $n$ such that $n\geq km$, there are
polynomials $T_1$ and $T_2$ in $S^{n}(B_p)$ such that $i_m(b_m)=
T_1/T_2$. In particular, such polynomials $T_1$ and $T_2$ exist in
$S^{km}(B_p)$.
\end{lemma}

\begin{proof}
The proof of Lemma~\ref{division} is similar to that of \cite[Proposition 1]{maclean2}. In \cite{volumes_chen}, Chen shows that if $B$ is approximable, then
there exist a constant d and another constant M such that
 \[ \dim(B_n) \sim M n^d.\]
In particular, there exist a constant $N$ and another constant $k$
such that for any $m_1, m_2> N'$ such that $m_2 > km_1$, we have
that
\[
\dim \left(B_{m_1+m_2}\right) \leq \tfrac{4}{3} \left(\dim
\left(B_{m_2}\right)\right).
\]
Pick now a $p$ and $n_0$ such that we have both of the
following:
\begin{enumerate}
 \item $p> N$, 
 \item $\dim(S^n(B_p)) \geq \frac{2}{3}{\dim (B_{np})}$ for all $n>n_0$. 
 \end{enumerate}
Now consider an element $b_m$ in $B_m$ for some $m>{\max}\{ N,
n_0p\}$. We assume furthermore that $m$ is divisible by $p$; \textit{i.e.}\ $m=
k'p$. Our aim is to give a bound on the poles of $i_m(b_m)$ which
depends linearly on $m$.

Choose $n$ such that $np \geq km$. We note that in particular $np \geq
N$.  We may assume that $i_m(b_m)= \frac{b_m}{b_m'}$ for some $b_m,
b_m'\in B_m$ since we have assumed that $B_1 \neq \{0\}$. Note that
\[\dim(b_m \cdot S^{n}(B_p))= \dim (S^n(B_p)) > \frac{2(\dim(B_{np}))}{3} > \frac{\dim(B_{np+m})}{2}.\]
Similarly,
 \[ \dim(b_m' \cdot S^{n}(B_p))= \dim (S^n(B_p)) > \frac{2(\dim(B_{np}))}{3}
 > \frac{\dim(B_{np+m})}{2},\]
 from which it follows that 
\[ b_m\cdot S^{n}(B_p))\cap b_m' \cdot  S^{n}(B_p))\neq \{0\}\]
 and hence 
\[ \frac{b_m}{b_m'}= \frac{T_1}{T_2}\]
 for some $T_1, T_2 \in S^n(B_p)$. This completes the proof of Lemma~\ref{division}.
\end{proof}
 
We can now complete the proof of Theorem~\ref{main}. Fix $n$ such that
$np=km$.  For $j=1, 2$ we set
 \[ (D(i_{np}(T_j)))= Z_j- P_j,\]
where $Z_j$ and $P_j$ are effective divisors that do not have any
common component; since $i_{np}(T_j)$ is a rational function, we have
in particular that $Z_j$ and $P_j$ are numerically equivalent. Note
that by the definition of $D_p$, we have that
 \[ P_j \leq (km)D_p,\]
and it follows that, numerically,
\[ Z_j\cdot H^{d-1} \leq (km) D_p\cdot H^{d-1}.\]
But now if we consider the poles divisor of $i_m(b_p)$, we have that
\[ P(i_m(b_m))\leq P_1+Z_2\]
and hence
\[P(i_m(b_m))\cdot H^{d-1} < 2(km)D_p\cdot H^{d-1}.\]
We know that there exists a $b_m$ such that $P(i_m(b_m))= D_m$ and
hence
\[\frac{D_m}{m} \cdot H^{d-1}\leq \frac{2km}{m} D_p\cdot
H^{d-1}\leq 2k D_p\cdot H^{d-1},\]
and hence this sequence is bounded since $k$ and $p$ were fixed.

It follows that the sequence $\frac{D_m}{m} \cdot H^{d-1}$ is bounded
and hence convergent, so that by Lemma~\ref{num_conv}, the sequence
$[D_m/m]$ is also convergent in $\NS(X)$. This completes the proof of Theorem~\ref{main}.\qed

%%%%%%%%%%%%%%%%%%%%%
% References
%%%%%%%%%%%%%%%%%%%%%

\end{document}